\documentclass[11pt,twoside]{article}
\usepackage{latexsym}
\usepackage{amssymb,amsbsy,amsmath,amsfonts,amssymb,amscd}
\usepackage{epsfig, graphicx}
\newcommand{\fer}[1]{(\ref{#1})}

\newcommand{\commentout}[1]{}
\newcommand{\R}{\mathbb{R}}

\newcommand {\e}  {\varepsilon}

\newcommand {\vp} {\varphi}

\newcommand {\Chi} {{\bf \raise 2pt \hbox{$\chi$}} }

\def\xb{{\bar x}}
\def\yb{{\bar y}}
\def\tb{{\bar t}}
\def\xe{{x_\e}}
\def\ye{{y_\e}}
\def\te{{t_\e}}
\def\ze{{z_\e}}
\newcommand {\f}   {\displaystyle\frac}
\newcommand {\p}   {\partial}
\newcommand{\dis}{\displaystyle}
\newcommand {\proof} {\noindent {\bf Proof}. }
\newcommand{\beq}{\begin{equation}}
\newcommand{\beqa}{\begin{eqnarray}}
\newcommand{\bea} {\begin{array}{rl}}
\newcommand{\beqan}{\begin{eqnarray*}}
\newcommand{\eeq}{\end{equation}}
\newcommand{\eeqa}{\end{eqnarray}}
\newcommand{\eeqan}{\end{eqnarray*}}
\newcommand{\eea} {\end{array}}
\newtheorem{theorem}{Theorem}[section]

\newtheorem{remark}[theorem]{Remark}

\newcommand{\qed}{{ \hfill
                       {\unskip\kern 6pt\penalty 500
                       \raise -2pt\hbox{\vrule\vbox to 6pt{\hrule width 6pt
                       \vfill\hrule}\vrule} \par}   }}
\title{\Large \bf Dirac concentrations in Lotka-Volterra parabolic PDEs}

\author{Beno\^ \i t Perthame\thanks{
Universit\'e Pierre et Marie Curie-Paris 6, UMR 7598 LJLL, BC187, 4, place Jussieu,  F-75252 Paris 5, and Institut Universitaire de France. Email: perthame@ann.jussieu.fr}
\and Guy Barles \thanks{Laboratoire de Math\'ematiques et Physique Th\'eorique, CNRS UMR 6083, F\'ed\'eration Denis Poisson,
Universit\'e Fran\c{c}ois Rabelais, Parc de Grandmont,
37200 Tours, France.
Email: barles@lmpt.univ-tours.fr
 }}

\date{\today}

\begin{document}
\maketitle
\pagestyle{plain}
%\tableofcontents
\pagenumbering{arabic}

\begin{abstract} 
We consider parabolic partial differential equations of Lotka-Volterra type, with a non-local nonlinear term. This models, at the population level, the darwinian evolution of a population; the Laplace term represents mutations and the nonlinear birth/death term represents competition leading to selection.

Once rescaled with a small diffusion, we prove that the solutions converge to a moving Dirac mass. The velocity and weights cannot be obtained by a simple expression, e.g., an ordinary differential equation. We show that they are given by a constrained Hamilton-Jacobi equation. This extends several earlier results to the parabolic case and to general nonlinearities. Technical new ingredients are a $BV$ estimate in time on the non-local nonlinearity, a characterization of the concentration point (in a monomorphic situation) and, surprisingly, some counter-examples showing that jumps on the Dirac locations are indeed possible. 
\end{abstract}

\\

AMS Class. No: 35B25, 35K57, 49L25, 92D15

%%%%%%%%%%%%%%%%%%%%%%%%%%%%%%%%%%%%%%%%%%%%
\section{Introduction}
\label{sec:intr}
%-------------------------------------------
%%%%%%%%%%%%%%%%%%%%%%%%%%%%%%%%%%%%%%%%%%%%

This paper is devoted to study the asymptotic behavior of solutions to parabolic Lotka-Volterra equations. They describe the dynamics of a population density $n(t,x)$ which expands (or decays) with a growth rate $R\big(x,I(t)\big)$ which changes sign. In the theory of adaptive evolution (\cite{metz,geritz,geritz2,MG,diekmann,CC}), the parameter $x$ denotes a physiological ``trait'' and $I(t)$ an environmental unknown shared by the total population and which is used as a nutrient. The population can use the nutrient $I(t)$ differently, depending upon the trait $x$, which makes that certain fitter traits should emerge, modifying the environment and thus allowing other traits to emerge. This mechanism uses mutations in the population that we model here by a mere diffusion; but more realistic integral kernels can be handled as well, \cite{DJMP,CCP}.  This elementary modeling leads to the equations  
\beq \left\{ \begin{array}{l}
\f{\p}{\p t} n_\e - \e \Delta n_\e =\f{n_\e}{\e} \; R\big(x,I_\e(t)\big), \qquad x\in \R^d, \; t \geq 0, 
\\ \\
n_\e(t=0)= n^0_\e \in L^1(\R^d) , \qquad n^0_\e\geq 0, 
\end{array} \right.
\label{eq:master}
\eeq
\beq 
I_\e(t)= \int_{\R^d} \psi(x) n_\e(t,x) dx. 
\label{eq:envir}
\eeq
The function $\psi$ is given and measures the 'predation' of individuals with trait $x$ on the environment. Notice that we have rescaled the problem in order to include the idea that mutations are small (or rare) thanks to the small parameter $\e>0$. Such models, together with related asymptotic, can be derived form individual based stochastic processes in the limit of large populations; we refer to \cite{champagnat, champagnatFM}. Also the problem without diffusion (a coupled infinite system of elementary differential equations) is interesting from the point of view of large time behavior; we expect that the dynamics concentrates on large times and several related results 
can be found in the literature, see \cite{DJMR, Pe} for instance.

Our purpose is to show that, under various assumptions on the rate $R(x,I)$, the population $n_\e$ concentrates as a Dirac mass (or a sum of Dirac masses), a mathematical way to express that well identified species emerge from the adaptive landscape defined by the rate $R$. The interesting feature being to describe the dynamics of these Dirac masses. Here, we perform this analysis in a completely rigorous manner. But the idea to analyze adaptive dynamics in those terms, and the formalism, goes back to \cite{DJMP}, where a well founded biological system (the chemostat) was studied. In  \cite{barlesP}, we performed a rigorous asymptotic analysis in the case of integral operators and with linear dependency on the environmental unknown $I(t)$. The case of a system, for a population with adults and juveniles, was studied in \cite{CCP}. 

Although related, the situation of reaction-diffusion systems, as they arise in combustion, is quite different. The simple Fisher-KPP equation may serve as a model, which amounts to use $R=1-n_\e$ in \fer{eq:master}, and is known to lead to the propagation of a front. This means that in the limit $\e \to 0$, the solution converges typically to either $0$ (uncolonized region) or to $1$ (fully colonized region), and the transition occurs on the ``front'', i.e., a hypersurface which dynamics can be described by the level set of a solution of a Hamilton-Jacobi Equation (see \cite{BES,BS1,BS2,ES,S1,freidlin}). Therefore, the limiting objects are, geometrically at least,  very different in the cases of reaction-diffusion equations and of Lotka-Volterra equations. This is the reason why in the later case a new type of equation, the {\em constrained Hamilton-Jacobi equation}, occurs to focus on the isolated points of the Dirac locations. The theory of viscosity solutions to Hamilton-Jacobi equations occurs naturally in our derivation because a phase occurs naturally; it is quite elaborate now and general introductory references are \cite{BCD, CIL,barlesbook, evans,FS}. But a large part of the present paper uses other general ideas that can be read without knowledge of this notion. 

The paper is organized as follows. We first state (section \ref{sec:amr})  simple and general results and in particular the convergence result to a Dirac mass for dimension one in ``monomorphic situations''. These results are completed by $BV$ bounds on $I_\e(t)$ (section \ref{sec:bv}) and by the asymptotic analysis through constrained Hamilton-Jacobi equations (section \ref{sec:aa}), in any dimension. Consequences are drawn in section \ref{sec:sp}, they imply the concentration of $n_\e$ as a sum of Dirac masses in a general multidimensional setting. Section \ref{sec:disc} is devoted to counterexamples to the continuity of $I(t)$ (the limit of $I_\e$) and of the dirac locations, and to ``smallness'' conditions implying continuity. Several further results are presented afterwards: the monomorphic situation is completed in section \ref{sec:mono}, the case of several environmental variables is treated in section \ref{sec:seu}.

%%%%%%%%%%%%%%%%%%%%%%%%%%%%%%%%%%%%%%%%%%%%
\section{Assumptions and main results}
\label{sec:amr}
%-------------------------------------------
%%%%%%%%%%%%%%%%%%%%%%%%%%%%%%%%%%%%%%%%%%%%

In this section, we give assumptions for the coefficients and data arising in \fer{eq:master}--\fer{eq:envir} and state the main convergence  result in a particularly simple case. More general results follow from the analysis and proof we perform later; these results, which are more technical, are stated along the paper.
\\

We assume that there are two constants $\psi_m$, $\psi_M$ such that 
\beq 
0<\psi_m \leq \psi \leq \psi_M < \infty,  \qquad \psi \in W^{2,\infty}(\R^d) . 
\label{aspsi}
\eeq
The quantity $R$ is called the {\em invasion exponent} because it describes the ability of the individuals of trait $x$ to invade the population with environmental state $I(t)$. 
We assume that there are two constants $0< I_m \leq  I_M < \infty$ such that
\beq
\min_{x\in\R^d} R(x, I_m)= 0, \qquad \max_{x\in\R^d} R(x, I_M) = 0,
\label{asr1}
\eeq
and there exists a constant $K>0$ such that, for any $x\in \R^d,$ $I \in \R$,
\beq 
-K \leq \f{\p R}{\p I} (x,I) < - K^{-1} <0,    \qquad \sup_{ I_m/2 \leq I \leq 2 I_M } \|R(\cdot, I)\|_{W^{2,\infty}(\R^d)} \leq K.  
\label{asr2}
\eeq
We will also use the assumption
\beq 
 n_\e^0 \in L^\infty(\R^d), \quad \nabla n_\e^0 \in L^1(\R^d) \quad \hbox{and}\quad I_m \leq   \int_{\R^d} \psi(x) \; n_\e^0(x) dx \leq I_M,
\label{asro}
\eeq
and the  notation
\beq 
\varrho_\e(t)=  \int_{\R^d} n_\e(t,x) dx .
\label{defro}
\eeq
Notice that the assumption \fer{aspsi} and the bound $I_m\leq \int_{\R^d} \psi(x) \; n_\e(t,x) \; dx   \leq I_M$ imply
$$
\f{I_m}{ \psi_M} =: \rho_m \leq \int_{\R^d}   n_\e(t,x) dx   \leq  \rho_M:=\f{  I_M}{\psi_m}.
$$

One can have in mind the particular but more intuitive example
\beq
\psi \equiv 1,   \qquad  R(x,I) = b(x) Q_1(I) -d(x)  Q_2(I), 
\label{mie}
\eeq
with $Q_i \in C^1(\R)$ for $i=1,2$, and 
\beq 
Q_1'(\cdot) <0, \qquad Q_2'(\cdot) >0, \qquad Q_i >0. 
\label{asq1}
\eeq
\beq 
b\geq b_m>0 , \quad  d\geq d_m >0 \text{ and } \; b, \; d \in W^{2,\infty}(\R^d). 
\label{as:lip}
\eeq

We recall that we have the following existence result and a priori bounds (here $C$ denotes various constants which maybe different from line to line)

%----------------------
\begin{theorem} With the assumptions \fer{aspsi}--\fer{asro} and $I_m - C \e^2\leq I_\e(0) \leq I_M+ C \e^2$, there is a unique solution $n_\e \in C\big(\R^+; L^1(\R^d)\big)$, 
to equation \fer{eq:master}--\fer{eq:envir}, and it satisfies, 
\beq
I_m -C \e^2 \leq I_\e(t) \leq I_M + C \e^2.
\label{eq:bound1}
\eeq
\label{th:basic}
\end{theorem} 
A proof of existence can be found, for instance, in \cite{DJMR} and we do not recall it. The uniform bound can also be also be found in \cite{barlesP, Pe} (according to \cite{DJMR}  the lower bound is unessential because they are indirect ways to guarantee non-extinction a posteriori, we keep it here for the sake of simplicity). Section \ref{sec:bv} gives another and stronger uniform bound in time, a uniform $BV$ bound on $I_\e(t)$. Here, we just indicate the derivation of the upper bound in \fer{eq:bound1}. We have 
$$
\begin{array}{rl}
\f d {dt}Ê\int_{\R^d} \psi(x) \; n_\e(t,x) \; dx &=\e  Ê\int_{\R^d}  n_\e(t,x)\; \Delta \psi \; dx + \f 1 \e \int_{\R^d}  n_\e(t,x) R\big(x,I_\e(t)\big)
\\[4mm]
& \leq  C_1 \e I_\e(t) + \f {1}{ \e} \; I_\e(t) \; \max_{x\in \R^d} \f{R\big(x,I_\e(t)\big) }{\psi(x)}
\end{array}
$$
and, from the assumption \fer{asr2}, the right hand side becomes negative as soon as $I_\e(t)$ overpasses $I_M+ \f{C_1 \psi_M}{K} \e^2$ and the result follows.

We can state a very simple version of our results in the simple case when dimension $d$ is equal to $1$ and when, typically, the function $R(x,I)$ is monotone in $x$.

%----------------------
\begin{theorem} [Dimension $d=1$] We assume \fer{aspsi}--\fer{asro}, the technical condition on $n^0_\e$ in Theorem \ref{conv-uni} below,  and 
\beq
\forall I_m < I < I_M \; \text{ there is a unique $X(I)\in \R$ such that } R\big(X(I), I\big)=0.
\label{as:mono}
\eeq
Then, the solution $n_\e(t)$ to equation \fer{eq:master}--\fer{eq:envir} converges in the weak sense of measures  (see also Remark \ref{rk:nlim} below)
\beq
n_{\e_k} (t) {\;}_{\overrightarrow{\; k \rightarrow \infty \; }}\; \varrho(t) \delta(x-\bar x(t)), 
\label{eq:monolim}
\eeq
and we have 
$$
\bar x(t)= X(I(t)), \qquad R\big(\bar x(t), I(t)\big)=0,
$$
and the pair $(\bar x(t), I(t))$ satisfies the constrained Hamilton-Jacobi equation given later on (see Section \ref{sec:aa}).

When $R(x,I)$ has the special form insuring uniqueness in Theorem \ref{conv-uni}, then the full family $n_\e$ converges.  
\label{th:mono}
\end{theorem} 

Such a population is called {\em monomorphic} because a single trait is represented asymptotically. This is the general situation with a single environmental variable $I(t)$ (this is called the Competitive Exclusion Principle, \cite{diekmann}). To go further, several environmental variables can be introduced (see \cite{DJMP,MG} and our results in Section \ref{sec:seu}).

%%%%%%%%%%%%%%%%%%%%%%%%%%%%%%%%%%%%%%%%%%%%
\section{BV estimates on $I_\e(t)$}
\label{sec:bv}
%-------------------------------------------
%%%%%%%%%%%%%%%%%%%%%%%%%%%%%%%%%%%%%%%%%%%%

As a first step in our analysis of the limit $\e \to 0$ in \fer{eq:master}, we prove strong convergence of $I_\e(t)$. Therefore, we complete result of Theorem \ref{th:basic} by the:

%-----------------------------
\begin{theorem}\label{BVe} With the assumptions \fer{aspsi}--\fer{asro}, we have additionally to the uniform $L^1$ bound \fer{eq:bound1}, the local uniform BV and sub-Lipschitz bounds  
$$
\f d {dt} I_\e(t) \geq - \e C +e^{-K_2 t/\e}    \int \psi(x) n_\e^0(x) \f{R\left(x, I_\e^0\right)}{\e} \; ,
$$
$$
\f d {dt} \varrho_\e (t) \geq - Ct  +  \int (1+ \psi(x)) n_\e^0(x) \f{R(x, I_\e^0)}{\e}\; ,
$$
where $C$ and $K_2$ are positive constants. 

Consequently the extracted limits satisfy that $I(t)$ is nondecreasing as soon as there exist a constant $C$ independent of $\e$ such that 
\beq
  \int \psi(x) n_\e^0(x) \f{R\left(x, I_\e^0\right)}{\e} 
 \geq - C e^{o(1)/\e } .
\label{as4}
\eeq
Then, we also have, for all $T>0$, 
\beq
\int_0^T \int_{\R^d} n_\e (t,x) R\big(x,I_\e(t)\big)^2 dx \; dt \leq C \; \e \; (1+ \e T ).
\label{eq:rsquare} 
\eeq
\end{theorem} 
%-----------------------------------------

\begin{remark} 
Our condition \fer{as4} on the initial data leaves place for a possible initial layer. In order to avoid it, one may choose {\em well-prepared initial data} that satisfy $R\left(x, I^0_\e\right)=0$ (for instance by tuning the total mass of $n^0_\e$). When it is not fulfilled one might observe, e.g.  in numerical simulations, a fast variation of $I_\e(t)$ for $t\approx 0$. 
\end{remark}
%----------------
\begin{remark} From the bound  \fer{eq:rsquare}, we can deduce that the weak limit $n$ in $L^\infty \big(\R^+; M^1(\R^d) \big)$ of $n_{\e_k}$ is supported (a.e. in $t$) only at points $x$ such that $R\big(x, I(t)\big)$. We give a more precise result below. 
\end{remark}
%-----------------
\begin{remark} Notice also that one can prove following the same lines that the family $\int_{\R^d} n_\e(t,x) \Psi(x) dx$ is bounded in $BV_{\rm loc}$ for any $W^{2,\infty}$ test-function $\Psi$. The extracted limits $n_{\e_k}$ therefore converge a.e. in time and not merely in 
$w\text{--}L^\infty \big(\R^+; w\text{--}M^1(\R^d)\big)$.
\label{rk:nlim}
\end{remark}
%-------------------
\begin{remark} The interested reader can check that the proof below also extends previous results of \cite{Pe, DJMR} for the continuous differential system 
$$
\f{\p}{\p t} n  =n \; R\big(x,I(t)\big), \qquad I(t)= \int_\R n(t,x) dx.
$$
When $R$ is monotonic in $I$ and assuming there is a unique $\bar x$ such that $0= R(\bar x , I_M)= \max_x R(x, I_M)$, then as $t \to \infty$, we have $n(t)\rightharpoonup \rho_M \delta(x-\bar x)$ with $I_M=\rho_M \psi(\bar x)$.
\end{remark}

\proof We begin with the proof for $I_\e$, then we show \fer{eq:rsquare} and finally indicate the variants for the proving the result on  $\varrho_\e$.

Concerning $I_\e$, we multiply the equation by $\psi$ and integrate over $\R^d$, this yields
\beq
\f d {dt}  I_\e(t) = \e  \int  n_\e(t,x) \Delta \psi(x) + {\cal J}_\e(t) ,
\label{eq:dtro}
\eeq
with $ {\cal J}_\e(t)$ defined by
$$
 {\cal J}_\e(t) =\int \f{n_\e}{\e} \;  \psi(x) \; R\big(x,I_\e(t)\big) .
$$
The integration by parts is justified because $\psi \in W^{2,\infty}$ and, for $\e$ fixed, both  $n_\e$ and $\nabla n_\e$ belong to $L^1$ by easy a priori manipulations.
 
The result relies on an estimate on $ {\cal J}_\e(t)$. In the same way, we have
\beq
 \bea
\f d {dt}  {\cal J}_\e(t) =&   \int n_\e \Delta \big[ \psi(x) R\big(x,I_\e(t)\big) \big] dx+\int \f{n_\e}{\e^2} \psi(x) R\big(x,I_\e(t)\big)^2 dx
\\ \\
& + \int \f{n_\e}{\e} \psi(x) \f{\p}{\p I}R\big(x,I_\e(t)\big)  dx \; \f d {dt}  I_\e(t).
\eea 
\label{derJ}
\eeq
Now we use \fer{eq:dtro} to recover $ {\cal J}_\e(t)$ from  $\f d {dt}  I_\e(t)$ in the last term. And we notice the following properties~: by (\ref{asr2}) and (\ref{eq:bound1}), we have
$$ \int n_\e \Delta \big[ \psi(x) R\big(x,I_\e(t)\big)  \big]
 -  \int  n_\e(t,x) \Delta \psi(x) \,  \int n_\e(t,x) \psi(x)  \f{\p}{\p I}R\big(x,I_\e(t)\big) =O(1) 
\geq -K_1,
$$
Therefore, using again  (\ref{asr2}), 
\beq
\f d {dt}  {\cal J}_\e(t) = O(1) +\int \f{n_\e}{\e^2} \psi(x) R\big(x,I_\e(t)\big)^2 dx + \int \f{n_\e}{\e} \psi(x) \f{\p}{\p I}R\big(x,I_\e(t)\big)  dx \; {\cal J}_\e(t) .
\label{eq:calj}
\eeq
But  we have 
$$ 
\int n_\e(t,x) \psi(x) \f{\p}{\p I}R\big(x,I_\e(t)\big) \leq -K_2 <0\; .
$$
Using also that the second term in the right-hand side of the above equality is positive, we obtain
$$
\f d {dt} \big(  {\cal J}_\e(t)\big)_- \leq K_1  -   \f{K_2}{\e} \big(  {\cal J}_\e(t)\big)_-.
$$
From this differential inequality we find 
\beq
\big(  {\cal J}_\e(t)\big)_- \leq \e \f{K_1}{K_2}+ \big(  {\cal J}_\e(0)\big)_-  e^{-K_2 t /\e}.
\label{eq:jemimis}
\eeq 
The local $BV$ bound  on $I_\e(t)= \int \psi(x) n_\e(t,x)dx$, as well as its monotonicity in the limit,  follows by inserting this inequality in  \fer{eq:dtro}. 
\\

In order to prove \fer{eq:rsquare}, we ague as follows. We have from \fer{eq:dtro}
\beq
\int_0^T{\cal J}_\e(t) = I_\e(T) -I_\e(0) +\e \; O(T)= C +\e \; CT.
\label{eq:jebound}
\eeq 
On the other hand, from \fer{eq:calj}, we also have, for some $K(t)>0$ bounded from above and from below away from $0$
$$
\f d {dt}  {\cal J}_\e(t) = C +\int_{\R^d} \f{n_\e}{\e^2} \psi(x) R\big(x,I_\e(t)\big)^2 dx - \f{K(t)}{\e} \; {\cal J}_\e(t) .
$$
After integration, this gives
$$
{\cal J}_\e(t) = {\cal J}_\e(0) e^{-\int_0^t  \f{K}{\e}}   + \dis \int_0^t e^{-\int_s^t  \f{K}{\e}}\left( C+ \int_{\R^d} \f{n_\e(s,x)}{\e^2} \psi(x) R\left(x,I_\e(s)\right)^2 dx\right) ds,
$$
and thus, using assumption \fer{as4} and with $ \underline K =\min K(t)>0$, we arrive at
$$
\int_0^T {\cal J}_\e(t) dt \geq -C+ \dis \int_{s=0}^T\left( C+ \int_{\R^d} \f{n_\e(s,x)}{\e^2} \psi(x) R\left(x,I_\e(s)\right)^2 dx\right) \; \f{\e}{\underline K} \; ds.
$$
According to \fer{eq:jebound}, we deduce that 
$$
\dis \int_{s=0}^T \int_{\R^d} \f{n_\e(s,x)}{\e} \psi(x) R\left(x,I_\e(s)\right)^2 dx\; ds \leq C + C \e \; T,
$$ 
and the result \fer{eq:rsquare} follows. 
\\

The result  on $\varrho_\e(t)=\int n_\e(t,x)dx$  follows similar lines. We  write
$$
\f d {dt} \varrho_\e(t)= \f 1 \e \int n_\e(t,x) R\big(x,I_\e(t)\big) dx :=  {\cal K}_\e(t) ,
$$
and 
\beq
 \bea
\f d {dt} {\cal K}_\e(t)& =   \int n_\e \Delta R\big(x,I_\e(t)\big) +\int \f{n_\e}{\e^2} R\big(x,I_\e(t)\big)^2 + \int \f{n_\e}{\e} \f{\p}{\p I}R\big(x,I_\e(t)\big)  \; \f d {dt}  I_\e(t)
\\ \\
& \geq -K_3+\f{ K_2}{\e}  {\cal J}_\e(t)
\\ \\
& \geq - K_1-K_3+ \f{ K_2}{\e} \big(  {\cal J}_\e(0)\big)_-  e^{-K_2 t /\e}
\eea 
\label{derK}
\eeq
(after using the inequality \fer{eq:jemimis} for $\big(  {\cal J}_\e(t)\big)_-$).
This implies
$$
 {\cal K}_\e(t) \geq  {\cal K}_\e(0) - (K_1 + K_3) t - \big(  {\cal J}_\e(0)\big)_-  .
$$
The result on $\f d {dt} \varrho_\e(t)$ follows. 
\qed

%%%%%%%%%%%%%%%%%%%%%%%%%%%%%%%%%%%%%%%%%%%%
\section{Constrained Hamilton-Jacobi equation, uniqueness}
\label{sec:aa}
%-------------------------------------------
%%%%%%%%%%%%%%%%%%%%%%%%%%%%%%%%%%%%%%%%%%%%

As already mentioned earlier, in the limit $\e \to 0$, the solution $n_\e$ to \fer{eq:master} converges weakly to a measure $n\in L^\infty\big(\R^+; L^1(\R^d)\big)$ (see see Remark \ref{rk:nlim}).  In this section we give a general theory for describing properties of $n$. The statements of Theorem \ref{th:mono} follow from the present analysis (and consequences in next section. 
\\

We expect that $n_\e$  concentrates as Dirac masses
$$
n_\e (t,x) \rightharpoonup n(t,x) = \sum_i \varrho_i(t) \delta \big(x- x_i(t)\big).
$$
This weak limit can be described more accurately through the phase function $\vp_\e$ defined by 
\beq
n_\e(t,x)= e^{\vp_\e(t,x)/\e},
\label{eq:phase} 
\eeq
just as the Dirac mass at $0$ is well approximated by the gaussian $\f{1}{\sqrt{2 \pi \e}} e^{-|x|^2/(2 \e)}$.  The description of $\vp_\e$, and its limit, gives information on the measure $n$.

Firstly, we obtain the following equation, equivalent to \fer{eq:master}, 
\beq 
\left\{ \begin{array}{l}
\f{\p}{\p t} \vp_\e(t,x) = |\nabla \vp_\e|^2+ R\big(x,I_\e(t)\big) +\e \Delta \vp_\e, 
\\ [4mm]
\vp_\e(t=0,x) = \vp^0_\e(x) := \e \ln n^0_\e.
\end{array} \right.
\label{eq:chjd}
\eeq
Following \cite{DJMP, barlesP}, in the limit, we obtain a viscosity solution to the constrained Hamilton-Jacobi equation
\beq
\left\{ \begin{array}{l}
\f{\p}{\p t} \vp(t,x) = |\nabla \vp|^2+ R\big(x,I(t)\big), 
\\  [4mm]
\dis \max_{x \in \R} \; \vp(t,x)= 0, \quad \forall t >0,
\\  [4mm]
\vp(t=0,x) = \vp^0(x).
\end{array} \right.
\label{eq:chj}
\eeq

%-------------------
\begin{theorem}\label{conv-uni} Assume \fer{aspsi}--\fer{asro} and $(\vp_\e^0)_\e$ is a sequence of uniformly bounded function in $W^{1, \infty}$ which converges uniformly to $\vp^0$. Then, after extraction of a subsequence, $(\vp_\e^0)_\e$ converges locally uniformly to a Lipschitz continuous viscosity solution $\vp$ to \fer{eq:chj}. In particular, a.e. in $t$, $supp \;n(t, \cdot) \subset \{ \vp(t,\cdot)=0\}$.

This solution is unique when assumptions \fer{as:lip} and
$$
R(x,I)= b(x)-d(x) Q(I), \qquad \text{with $Q(I)>0$ increasing},
$$
or 
$$
R(x,I)= b(x)Q(I)-d(x), \qquad \text{with $Q(I)>0$ decreasing},
$$
\end{theorem}
%----------------
Because we expect that the set $ \{ \vp(t,\cdot)=0\}$ is made of isolated points, we indeed expect that $n$ is a sum of Dirac masses. In one dimension this follows rigorously from Theorem~\ref{th:lscx} below but heuristically maximum points of $\vp(t, \cdot)$ are indeed isolated in any dimension. 
\\

The result was proved in \cite{barlesP} in the case $R(x,I)= b(x)-d(x) Q(I)$ and for an integral operator instead of a diffusion. For completeness, we recall the proof and the main new ingredients. One of the difficulties being that coefficients, because of $I(t)$, are not discontinuous. This can be handled following well established arguments, see
\cite{barlesdisc, barlesP}. 
\\

\proof {\em Existence.}  In order to prove it, we show how to pass to the limit in (\ref{eq:chjd}). It is worth pointing out that we can do it in two ways, either by using the notion of viscosity solutions for equations with a $L^1$-dependence in time or, since $I_\e$ converges to an increasing function, by using the notion of viscosity solutions for equations with a discontinuous Hamiltonians. We do it in the second way which has the advantage to give a more precise result.

Since the functions $\vp_\e$ are equi-bounded and equi-Lipschitz continuous in space, they are also equi-H\"older continuous in time (See \cite{bbl} or \cite{genieysP} for proofs of this claim). Therefore the only difficulty is to pass to the limit (up to a subsequence) in the term $R(x,I_\e)$.

More precisely, in order to pass to the limit in this term, we have to compute the quantities (up to a subsequence)
$$ \limsup_{\displaystyle{\mathop{\scriptstyle{s\to
t}}_{ \e' \to 0}}}\, I_{\e'} (s) \quad \hbox{ and }\liminf_{\displaystyle{\mathop{\scriptstyle{s\to
t}}_{ \e' \to 0}}}\, I_{\e'} (s) \; .$$
Indeed, the $x$-dependence in $R(x,I_\e)$ does not cause any problem since the functions $x \mapsto R(x,I)$ are Lipschitz continuous in $x$ and, on the other hand, the function $R(x,I)$ is decreasing in $I$ for any $x$ which allows to reduce the computation of these limsup and liminf to $I_{\e'}$ and not $R(x,I_{\e'})$.

To do so, we use Theorem~\ref{BVe}~: if we set $\mu_\e := {1\!\!1}_{[\e^{1/2},T]}\left(\frac{dI_\e}{dt} + 2C\e\right)dt$, where, here and below, ${1\!\!1}_{[a,b]}$ denotes the indicator function of $[a,b]$, then the $\mu_\e$'s are bounded (positive) measures on $(0, T)$ for any $T>0$, if $\e$ is small enough. Therefore we can extract a subsequence such that $(\mu_{\e'})_{\e'}$ converges weakly to some measure $\mu$.

If $0< t_0 < T$ is a point such that $\mu(\{t_0\})=0$, then we can also assume that $I_{\e'} (t_0)$ converges to some $\alpha \in \R$ and we set
$$ I(t) = \alpha + \int_{(t_0,t)}\, d\mu(s)\quad \hbox{if  }t >t_0 \quad\hbox{and}\quad I(t) = \alpha - \int_{(t,t_0)}\, d\mu(s)\quad \hbox{if  }t \leq t_0\; .$$
Obviously, the function $I$ is an increasing function on $(0,T)$.

We are going to show that
$$ \limsup_{\displaystyle{\mathop{\scriptstyle{s\to
t}}_{ \e' \to 0}}}\, I_{\e'} (s) \leq I(t+) \quad \hbox{ and }\quad \liminf_{\displaystyle{\mathop{\scriptstyle{s\to
t}}_{ \e' \to 0}}}\, I_{\e'} (s) \geq I(t-) \; .$$
We do it for $t>t_0$, the other case being treated analogously.

For $s$ close to $t$, we have
$$I_{\e'} (s) = \alpha + \int_{(t_0,t)}\, d\mu_\e (\tau)\; .$$
If $\zeta: [0,T] \to \R $ is a continuous function such that $\zeta(\tau) \geq  {1\!\!1}_{[t_0,t]}(\tau)$ on $[0,T]$, we have
$$I_{\e'} (s) \leq  \alpha + \int_{(0,T)}\, \zeta(\tau) d\mu_\e (\tau)\; ,$$
and passing to the limit in this inequality, we obtain
$$ \limsup_{\displaystyle{\mathop{\scriptstyle{s\to
t}}_{ \e' \to 0}}}\, I_{\e'} (s) \leq \alpha + \int_{(0,T)}\, \zeta(\tau) d\mu (\tau)\; ,$$ and we conclude by applying Lebesgue's dominated convergence Theorem to a sequence $(\zeta_k)_k$ which converges to $ {1\!\!1}_{[t_0,t]}$ in a suitable way : namely, we choose a sequence of functions such that $0\leq \zeta_k \leq 1$ on $[0,T]$,  $\zeta_k \equiv 1$ on $[t_0,t]$ for any $k$ and which converge pointwise to $0$ outside the interval $[t_0,t]$.

The proof for the liminf is the same approximating, this time, the indicator function from below. It is worth pointing out that the property $\mu(\{t_0\})=0$ avoids here (and above) discussions on the behavior of the sequence at this point. We find that this liminf is larger than $ \alpha + \int_{(t_0,t)}\, \zeta(\tau) d\mu (\tau)$, which is exactly $I(t-)$.
\qed

\medskip

\proof {\em Uniqueness.} Consider for instance the second case. We consider again the function
$$
\Psi(t,x) = \vp(t,x) - b(x) \Sigma(t), \qquad \Sigma(t)=\int_0^t Q\big(I(s)\big) ds.
$$
It satisfies
$$
\f{\p}{\p t} \Psi(t,x) = -d(x) + \big|\nabla \big( \Psi +b(x)  \Sigma(t)\big)  \big|^2 .
$$
\\

On the one hand, for two different solutions with the same initial data, we define in this way two functions $\Psi_1$ and $\Psi_2$. Using the viscosity criteria, we have, at the point $x_0$ where the maximum is achieved 
\begin{eqnarray}
\f{d}{d t} \| \Psi_1 - \Psi_2(t) \|_\infty & = & \big|\nabla \big( \Psi_1(x_0) +b(x_0)  \Sigma_1(t)\big) \big|^2-\big|\nabla \big( \Psi_2(x_0) +b(x_0)  \Sigma_2(t) \big) \big|^2 \nonumber
\\  
&= & \big(\nabla\Psi_1(x_0) +\nabla b(x_0)  \Sigma_1(t) + \nabla\Psi_2(x_0) +\nabla b(x_0)  \Sigma_2(t) \big).\big(\nabla b(x_0)  \Sigma_1(t) - \nabla b(x_0)  \Sigma_2(t) \big) 
\nonumber 
\\  
&\leq & C |\nabla b(x_0)  \Sigma_1(t) - \nabla b(x_0)  \Sigma_2(t) | 
\nonumber 
\\  
&\leq & C \|Ê\nabla b\|_\infty \; |  \Sigma_1(t) -  \Sigma_2(t) | \label{dtineq}
\end{eqnarray}

On the other hand, we also have, considering the point $x_i$ where  $ \max_{x \in \R^d} \; \vp_i(t,x)$ is achieved
\begin{eqnarray}
0 & = &\dis \max_{x \in \R^d} \; \vp_1(t,x) - \dis \; \max_{x \in \R^d} \vp_2(t,x) \nonumber \\
& \leq &   \vp_1(t,x_1) - \vp_2(t,x_1)
\nonumber \\
& \leq & b(x_1)[\Sigma_1(t)- \Sigma_2(t)] +  \Psi_1(t,x_1) - \Psi_2(t,x_1)
\nonumber \\
& \leq & b(x_1)[\Sigma_1(t)- \Sigma_2(t)] +  \dis \max_{x \in \R^d} \; [\Psi_1(t,x) - \Psi_2(t,x)].
\label{ineqpsi1}
\end{eqnarray}

Changing the indices $1$ and $2$, one of the numbers $b(x_1)[\Sigma_1(t)- \Sigma_2(t)] $ or $b(x_2)[\Sigma_2(t)- \Sigma_1(t)] $ is  negative and thus we have, with the notation \fer{as:lip},  
\begin{equation}\label{ineqpsi2}
b_m \big| (\Sigma_1(t)- \Sigma_2 (t) \big| \leq  \dis \max_{x \in \R^d} \;  \big| \Psi_1(t,x) - \Psi_2(t,x)  \big|.
\end{equation}

Together, the two above inequalities yield
\beq
\f{d}{d t} \| (\Psi_1 - \Psi_2) (t) \|_\infty \leq C \| (\Psi_1 - \Psi_2) (t) \|_\infty. 
\label{eq:uniqgronw}
\eeq
Uniqueness follows.
\qed

\medskip

The uniqueness result can be obtained under a slightly more general assumption on $R$, namely:
\\[2mm]
for any functions $I_1, I_2 \in L^\infty (0,T)$ and $t\in (0,T)$, the function $x \mapsto \int_0^t \left[R\big(x, I_1(s)\big) - R\big(x, I_2(s)\big)\right]ds$ has a constant sign for all  $x \in \R^d$ and there exists a constant $K>0$ such that
\beq \label{hyp:Runi}
 \max_{x \in \R^d} \left | \int_0^t \left[D_x R\big(x, I_1(s)\big) - D_x R\big(x, I_2(s)\big)\right]ds\right | \leq \tilde K \min_{x \in \R^d}  \left | \int_0^t \left[R\big(x, I_1(s)\big) - R\big(x, I_2(s)\big)\right]ds\right |\; .
\eeq

The result can be stated as follows
%------------------------------------------
\begin{theorem}\label{uni-gen} Assume that $\vp^0 \in W^{1, \infty} (\R^d)$ and that \fer{aspsi}--\fer{asro}, \fer{hyp:Runi} and the above sign condition  hold. Then the problem \fer{eq:chj} has at most one Lipschitz continuous solution.
\end{theorem}

The two conditions appearing in assumption \fer{hyp:Runi}  are easily checkable in the cases we emphasize in Theorem~\ref{conv-uni} but (a priori) they are not satisfied for general $R$ of the form $R(x,I) = b(x) Q_1(I) -d(x)  Q_2(I)$, even if (\ref{asq1}), (\ref{as:lip}) hold true. Uniqueness for such general $R$ is still an open problem.

\medskip

\proof The proof follows readily the uniqueness proof of Theorem~\ref{conv-uni}. We consider two solutions $(\vp_1, I_1)$ and $(\vp_2, I_2)$ and for $i=1,2$, we set 
$$
\Psi_i (t,x)= \vp_i (t,x)- \int_0^t R\big(x, I_i (s)\big) ds.
$$

The first part of the proof leading to \fer{dtineq} yields now
$$ \f{d}{d t} \| \Psi_1 - \Psi_2(t) \|_\infty \leq K  \max_{x \in \R^d} \left | \int_0^t \left[D_x R\big(x, I_1(s)\big) - D_x R\big(x, I_2(s)\big)\right]ds\right | , $$
and the second part (cf. \fer{ineqpsi1}--\fer{ineqpsi2}) provides us with the inequalities
$$ 0 \leq  \int_0^t \left[ R\big(x_1, I_1 (s)\big)- R\big(x_1, I_2 (s)\big) \right]ds +  \dis \max_{x \in \R^d} \; [\Psi_1(t,x) - \Psi_2(t,x)],
$$
$$ 0 \leq  \int_0^t \left[R\big(x_2, I_2 (s)\big)- R\big(x_2, I_1 (s)\big)\right]ds +  \dis \max_{x \in \R^d} \; [\Psi_2(t,x) - \Psi_1(t,x)].
$$
and then
$$  \min_{x \in \R^d}  \left |  \int_0^t \left[R\big(x, I_1 (s)\big)- R\big(x, I_2 (s)\big)\right]ds\right |  \leq \dis \max_{x \in \R^d} \; | \Psi_1(t,x) - \Psi_2(t,x)|
$$
(this is where it is important that $x \mapsto \int_0^t \left[R\big(x, I_1(s)) - R(x, I_2(s)\big)\right]ds$ does not change sign).
Putting together the two inequalities for $\f{d}{d t} \| \Psi_1 - \Psi_2(t) \|_\infty$ and $\| \Psi_1 - \Psi_2(t) \|_\infty$ and using assumption \fer{hyp:Runi} leads to \fer{eq:uniqgronw} and uniqueness follows.
\qed

\medskip

\begin{remark} In \cite{barlesP}, the uniqueness proof is given in full details for merely continuous solutions. Lipschitz continuity of $\vp$ is used here to simplify the arguments (an write them a.e.) but everything can be understood in viscosity sense.  
\end{remark}

%%%%%%%%%%%%%%%%%%%%%%%%%%%%%%%%%%%%%%%%%%%%
\section{Structure properties}
\label{sec:sp}
%-------------------------------------------
%%%%%%%%%%%%%%%%%%%%%%%%%%%%%%%%%%%%%%%%%%%%
%%%%%%%%%%%%%%%%%%%%%%%%%%%%%%%%%%%%%%%%%%%%

We now derive consequences of the constrained Hamilton-Jacobi equation in terms the concentration points.
%--------------------------------------------
\begin{theorem}\label{propx} Assume \fer{asr2}. For any $\vp^0 \in W^{1, \infty}(\R^d)$, the solution to \fer{eq:chj} is semi-convex in $x$ for $t>0$, i.e. for any $t>0$, there exist a $C_{sc}(t)$ such that, for any $\xi \in \R^d$, $|\xi| =1$, we have 
$$
 \f{\p^2}{\p \xi^2} \vp_\e \geq -C_{sc}(t), \qquad i=1,...,d.
$$
 Consequently, at a maximum point $\bar x(t)$ of $\vp(t, \cdot)$ in the variable $x$, $\vp(t, \cdot)$ is differentiable in $x$ (but maybe not in $t$) and we have 
$$
 \nabla \vp\big(t,\bar x(t)\big)=0,
$$
additionally for all Lebesgue points of $I(t)$ we have
\beq
R\big(\bar x(t),I(t)\big) =0 .
\label{eq:iv}
\eeq
\label{th:lscx}
\end{theorem}

The following proof does not use the $BV$ property of $I(t)$, (only the strong limit $I_\e \to I$ is used to pass to the limit in the Hamilton-Jacobi equation), and thus might be useful for extensions of the model where $BV$ estimates are not proved  as in section \ref{sec:seu}. 
\\

\proof 
{\em (Proof of semi-convexity)} We justify the semi-convexity directly on the diffusive equation and use the notation $\vp_\xi:= \f{\p \vp_\e}{\p \xi} $, 
$\vp_{\xi\xi}:= \f{\p^2 \vp_\e}{\p \xi^2} $. We have 
$$
\f{\p}{\p t}Ê\vp_\xi = 2 \nabla \vp_\e \cdot \nabla \vp_\xi + R_\xi(x,I_\e(t)) + \e \Delta \vp_\xi, 
$$
\beq
\f{\p}{\p t}Ê\vp_{\xi\xi} = 2 \nabla \vp_\e \cdot \nabla \vp_{\xi\xi} +2 |\nabla \vp_\xi|^2+ R_{\xi\xi}(x,I_\e(t)) + \e \Delta \vp_{\xi\xi}.
\label{eq:secondder}
\eeq
But $|\nabla \vp_\xi| \geq |\vp_{\xi\xi}|$ because $\vp_{\xi\xi}= \nabla \vp_\xi\cdot \xi$, therefore the function $w:=\vp_{\xi\xi}$ satisfies
$$
\f{\p}{\p t}Êw \geq 2 \nabla \vp_\e \cdot \nabla w + 2 w^2  -  \bar R+ \e \Delta w
$$
with $\bar R$ an upper bound on $R_{\xi\xi}(x, I_\e(t))$. The semi-convexity result (in fact, the semi-convex regularizing effect) follows from suitable approximation arguments (for the initial data) and the comparison with the (absolute) subsolution given by the solution of the O.D.E. $\dot y = 2 y^2  - \bar R$, $y(0) = -\infty$.
\\ 
\\
{\em (Proof of $\nabla \vp\big(t,\bar x(t)\big) = 0$)} We prove this equality and a related preliminary result. As we have shown above, the function $\vp$ is semi-convex and, by classical properties of semi-convex functions, $\vp$ is differentiable at maximum points. 
Therefore
$$
\nabla \vp\big(t, \bar x(t)\big)= 0. 
$$
Moreover, it is standard that for any sequence $(t_k, x_k)$ of $x$-differentiability points of $\vp$ which converges to $(t, \bar x(t))$, we also have 
$$ \nabla \vp (t_k, x_k) \to 0. 
$$
From this property, we deduce that, for $h, r >0$, $h, r \to 0$, 
\beq
 \frac{1}{rh}\int_t^{t+h}\int_{ \bar x(t)-r}^{ \bar x(t)+r} |\nabla \vp(s,y)|^2 dsdy \ , \ \frac{1}{rh}\int_{t-h}^{t}\int_{ \bar x(t)-r}^{ \bar x(t)+r} |\nabla \vp(s,y)|^2 dsdy \to 0\; .
 \label{psc}
 \eeq
Indeed, considering the first limit for instance, a straightforward change of variable yields
$$
\int_0^{1}\int_{ -1}^{ 1} |\nabla \vp\big(t + h \tau , x(t)+ r e\big)|^2 d\tau de\; ,
$$
and it suffices to apply Lebesgue's dominated convergence Theorem
to this integral, using the above property together with the (local) Lipschitz continuity of $\vp$.
\\
\\
{\em (Proof of $R\big(t,\bar x(t)\big)\leq 0$)} 
We come back to the equation which holds almost everywhere, and we first integrate it on rectangles $(t, t+h) \times ( \bar x(t)-r, \bar x(t)+r) $. We obtain
$$
\int_{ \bar x(t)-r}^{ \bar x(t)+r} (\vp(t+h,y) - \vp(t,y) dsdy \int_t^{t+h}\int_{ \bar x(t)-r}^{ \bar x(t)+r} R \big(y,I(s)\big) dsdy +
\int_t^{t+h}\int_{ \bar x(t)-r}^{ \bar x(t)+r} |\nabla \vp(s,y)|^2 dsdy\; .
$$
But, by the semi-convexity of $\vp$ in $x$, we have, for $y \in ( \bar x(t)-r, \bar x(t)+r)$
$$ 0 \geq \vp(t,y) \geq \vp(t,\bar x(t)) - \f{C_{sc}}{2} |y-\bar x(t)|^2 = O(r^2)\; ,$$
while we have also $\vp(t+h,y) \leq 0$.
Using these two properties in the above equality, we deduce
$$  
  \frac{1}{rh}\int_t^{t+h}\int_{ \bar x(t)-r}^{ \bar x(t)+r} R \big(y,I(s)\big) dsdy +
 \frac{1}{rh} \int_t^{t+h}\int_{ \bar x(t)-r}^{ \bar x(t)+r} |\nabla \vp(s,y)|^2 dsdy \leq  \frac{1}{rh}O(r^2)\; .
$$
Therefore, we obtain
$$
\frac{1}{rh}\int_t^{t+h}\int_{ \bar x(t)-r}^{ \bar x(t)+r} R \big(y,I(s)\big) dsdy \leq  \frac{1}{rh}O(r^2)
$$
and letting $r,h$ tend to zero with $r \ll h$, we conclude that at any Lebesgue point of $I$ we have
$$ R\big(\bar x(t), I(t)\big)\leq 0\; .
$$
\\
{\em (Proof of $R\big(t,\bar x(t)\big)\geq 0$)} 
To obtain the opposite inequality, we integrate on the rectangle $(t-h, t) \times ( \bar x(t)-r, \bar x(t)+r) $. The left-hand side is now
$$ 
\int_{ \bar x(t)-r}^{ \bar x(t)+r} (\vp(t,y) - \vp(t-h,y) dsdy\geq \int_{ \bar x(t)-r}^{ \bar x(t)+r} \vp(t,y)  dsdy 
$$
and, by the afore mentioned argument, it is larger than $O(r^2)$ and we are lead to
$$ 
 \frac{1}{rh}\int_{t-h}^t\int_{ \bar x(t)-r}^{ \bar x(t)+r} R \big(y,I(s)\big) dsdy +
 \frac{1}{rh} \int_{t-h}^t\int_{ \bar x(t)-r}^{ \bar x(t)+r} |\nabla \vp(s,y)|^2 dsdy \geq  \f{O( r )}{h} \; .
$$
Again we let $r,h$ tend to zero with $r \ll h$ and using  (\ref{psc}), we conclude that, at any Lebesgue point of $I$, we have
$$
 Q\big(t,\bar x(t)\big)\geq 0\; .
$$
\qed

\bigskip

Let us recall that because $I(t)$ is $BV$, it is continuous away from a countable set of discontinuity points. Hence, we can  give another proof of the relation \fer{eq:iv}. Let  $t_0$ a continuity point 
of $I(t)$ and let  $\bar x(t_0)$ be a maximum point of $\vp(t_0)$. 

First, we use the viscosity subsolution criteria at such a point $(t_0, \bar x(t_0))$, testing against the test function $0$. We find
\beq
0 \leq R\big( \bar x(t_0),I(t_0)\big) +0.
\label{eq:iv1}
\eeq
To get the other inequality , we integrate in time the equation \fer{eq:chj} on $(t_0,t_0+h)$ ($h>0$) at the point $x=\bar x(t_0)$ and find
$$
0 \geq \f{1}{h} \vp\big(t_0+h, \bar x(t_0)\big) \geq  \f{1}{h} \int_{s=0}^h R\big( \bar x(t_0), I(t_0+s)\big) ds.
$$
Because $t_0$  is a continuity point of $I(t_0)$, we find
\beq
0 \geq \liminf_{h \to 0^+} \f{ \vp\big(t_0+h, \bar x(t_0)\big)- \vp\big(t_0,\bar x(t_0)\big)}{h} \geq  R\big( \bar x(t_0),I(t_0)\big) .
\label{eq:iv2}
\eeq
Clearly, \fer{eq:iv1} and \fer{eq:iv2} give the relation \fer{eq:iv}.
\qed

%%%%%%%%%%%%%%%%%%%%%%%%%%%%%%%%%%%%%%%%%%%%
\section{Continuity and discontinuity of $I(t)$ and $\bar x (t)$}
\label{sec:disc}
%-------------------------------------------
%%%%%%%%%%%%%%%%%%%%%%%%%%%%%%%%%%%%%%%%%%%%
%%%%%%%%%%%%%%%%%%%%%%%%%%%%%%%%%%%%%%%%%%%%

So far our results have shown that the Lagrange multiplier $I(t)$ is $BV$ and that the concentration point $\bar x(t)$ is unique in one dimension. It is natural to ask whether these quantities are continuous. This question is particularly relevant because it is possible to derive formally a differential equation on $\bar x(t)$, see \cite{DJMP, barlesP}
$$
\f{d}{dt} \bar x(t) = \left( -D^2 \vp(t,x) \right)^{-1} \nabla_x R\big(\bar x(t), I(t)\big). 
$$

In this section, we answer negatively to this question (in general) based on analytical and numerical  examples. We also give a smallness assumption implying continuity. 

%-------------------------------------------------------------
\subsection{An analytical counter-example}
%-------------------------------------------------------------

Our analytical counterexample follows earlier ideas developed for combustion models, see \cite{freidlin, S1} for instance. 

We choose in this section $\psi \equiv 1$, in other words $I(t)=  \varrho(t)$. We first consider the simple equation
\beq \label{shj}
\vp_t = x -  \varrho(t) +|\vp_x|^2 \quad\hbox{in  } (0, +\infty)  \times  \R, 
\eeq
together with the initial data
\beq \label{id:shj}
\vp (x,0) = \vp_0 (x) \quad\hbox{in  } \R \; ,
\eeq
and the constraint
\beq \label{cons:shj}
\max_{x \in \R} \, \vp(t,x) = 0\quad\hbox{for any  }t\; .
\eeq
In a first step, we start by solving the initial value problem for $\Psi = \vp - \int_0^t  \varrho(s) ds$  with the initial data $  \vp_0 (x) = - x^2$. We look for a solution of the form
$$ a(t)(x- b(t))^2 + c(t)\; .
$$
Easy computations yield 
$$ \left\{ \begin{array}{ll}
 a'(t) = 4a^2 (t), & \quad \quad a(0) = -1\; ,\\ \\
 -2a(t) b'(t) = 1, & \quad \quad b(0) = 0\; ,
\\ \\
c'(t) = b(t), & \quad \quad c(0) = 0\; .
\end{array} \right.
$$
Therefore, we have the explicit coefficients and solution 
$$ a(t) = - \frac{1}{1+4t}\; ,\quad  b(t) = \frac{t}{2}+t^2\; , \quad  c(t) =  \frac{t^2}{4}+ \frac{t^3}{3}\; ,
$$
$$ \Psi(t,x)= -\frac{1}{1+4t}\left(x- \frac{t}{2}-t^2\right)^2 + \frac{t^2}{4}+ \frac{t^3}{3}\; .
$$
From this, we deduce that the solution to (\ref{shj})--(\ref{cons:shj}) is
$$ \vp (t,x) = -\frac{1}{1+4t}\left(x- \frac{t}{2}-t^2\right)^2 \; , \qquad   \varrho(t) = \frac{t}{2}+t^2\; .
$$
This is consistent with the property $ \varrho(t) = \bar x(t) = b(t)$ of Theorem \ref{th:mono}.

In a second step, we remark that, in order to solve the problem on $\Psi = \vp - \int_0^t  \varrho(s) ds$ with initial data $-(x - \alpha)^2 -\delta$ ($\alpha >0$), we have just (i) to translate the solution obtained with $\alpha = 0$, (ii) to subtract $\delta-\alpha t$ to take into account the new right-hand side of the equation and the initial data. Therefore this solution is
$$ -\frac{1}{1+4t}\left(x- \alpha - \frac{t}{2}-t^2\right)^2 - \delta +\alpha t + \frac{t^2}{4}+ \frac{t^3}{3}\; .$$

Finally, we solve (\ref{shj})--(\ref{cons:shj}) with
$$ \vp_0 (x) = \max ( -x^2, -(x - \alpha)^2 -\delta)\quad\hbox{in  } \R \; .
$$
Since the Hamiltonian of (\ref{shj}) is concave, the maximum of two solutions is a solution : this is a consequence of the Barron-Jensen approach (\cite{BJ1,BJ2}, see also \cite{barlesdisc}). Thus the solution of this problem is
$$ \max \left(-\frac{1}{1+4t}\left(x- \frac{t}{2}-t^2\right)^2 , -\frac{1}{1+4t}\left(x-  \alpha - \frac{t}{2}-t^2\right)^2 - \delta +\alpha t\right) \quad \hbox{for  } t\leq {\bar t}:= \frac{\delta}{\alpha}\; ,$$
with $ \varrho (t) =  \frac{t}{2}+t^2$, while
$$ \max \left(-\frac{1}{1+4t}\left(x- \frac{t}{2}-t^2\right)^2  + \delta -\alpha t, -\frac{1}{1+4t}\left(x-  \alpha - \frac{t}{2} - t^2\right)^2 \right) \quad \hbox{for  } t > {\bar t}\; ,$$
with $ \varrho (t) =\alpha + \frac{t}{2}+t^2$.
\\

This shows the phenomena of discontinuity of $ \varrho(t)$ and of the Dirac concentration point $\bar x(t)$.

%%%%%%%%%%%%%%%%%%%%%%
\begin{figure}[ht!]
\begin{center}
\includegraphics[height=45mm,width=65mm]{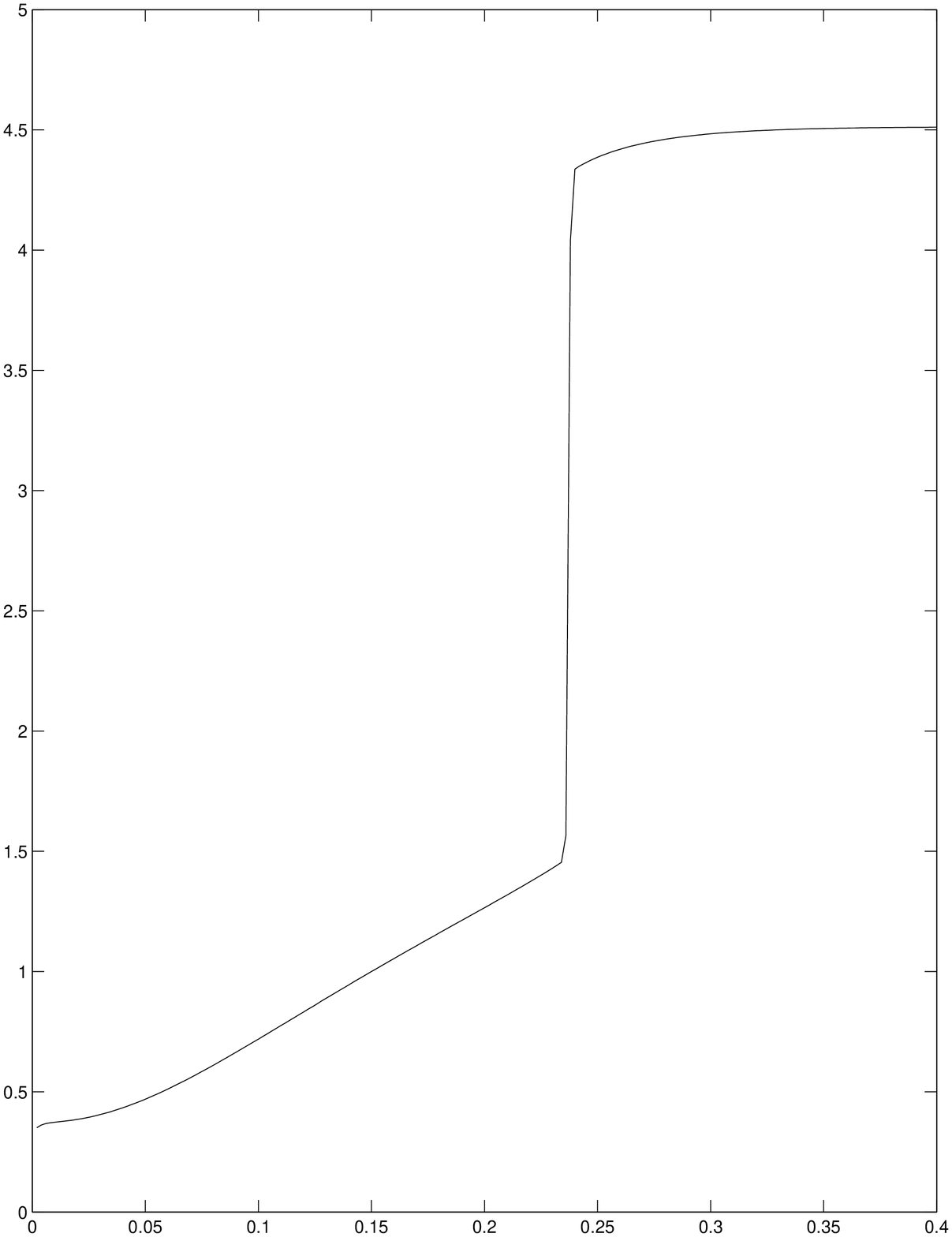} \quad \includegraphics[height=45mm,width=65mm]{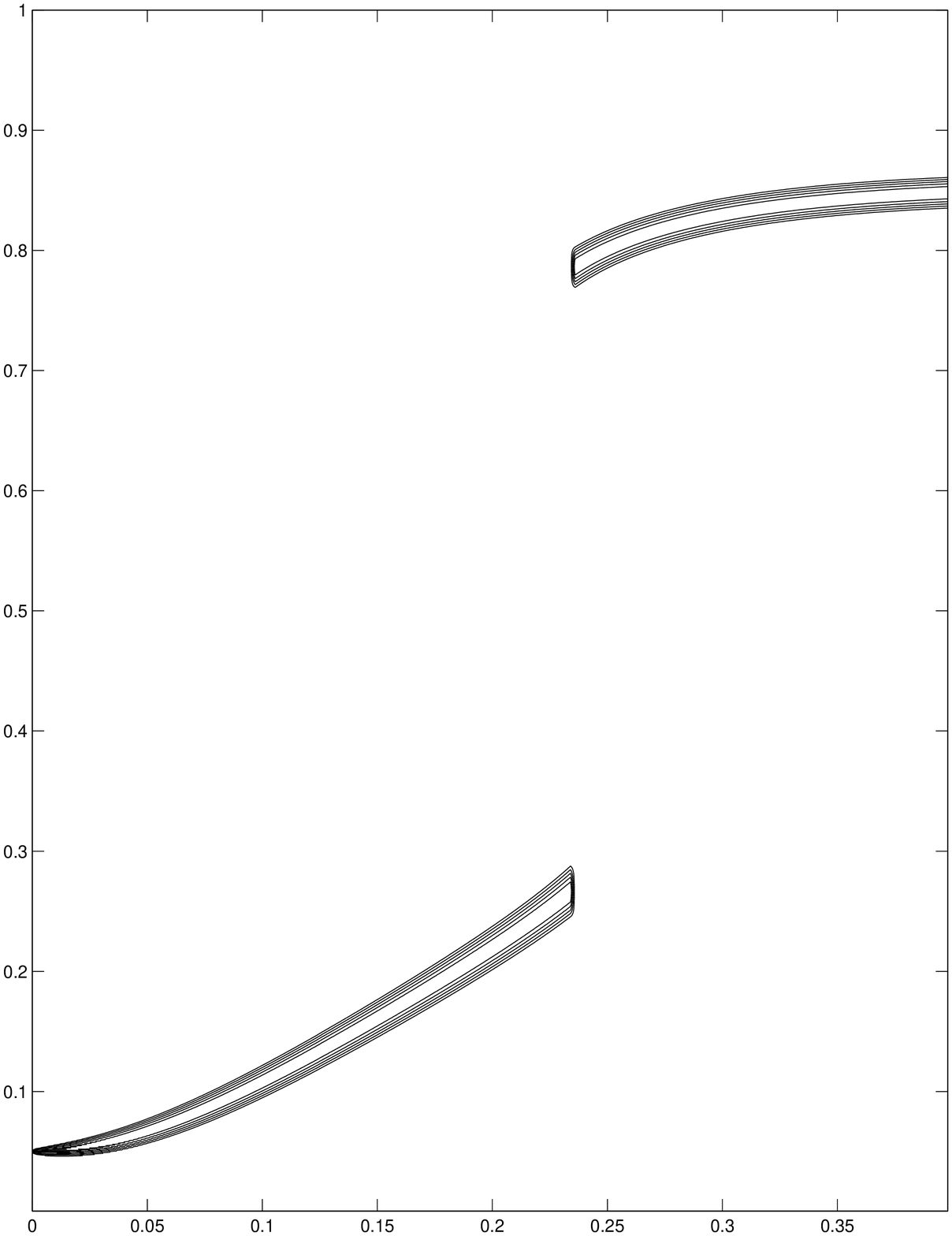}
\end{center}
\vspace{-.7cm}
\caption{
\label{fig:disc} 
{\footnotesize\sc 
A numerical example of discontinuous solution in example \fer{exR}. The abscissae is time. (Left)  The Lagrange multiplier $ \varrho(t)$. (Right) The concentration point $\bar x(t)$ (in fact isovalues of the density $n_\e(t,x)$).
}
        }
\end{figure}
%%%%%%%%%%%%%%%%%%%%%%
\begin{figure}[ht!]
\begin{center}
\includegraphics[height=45mm,width=50mm]{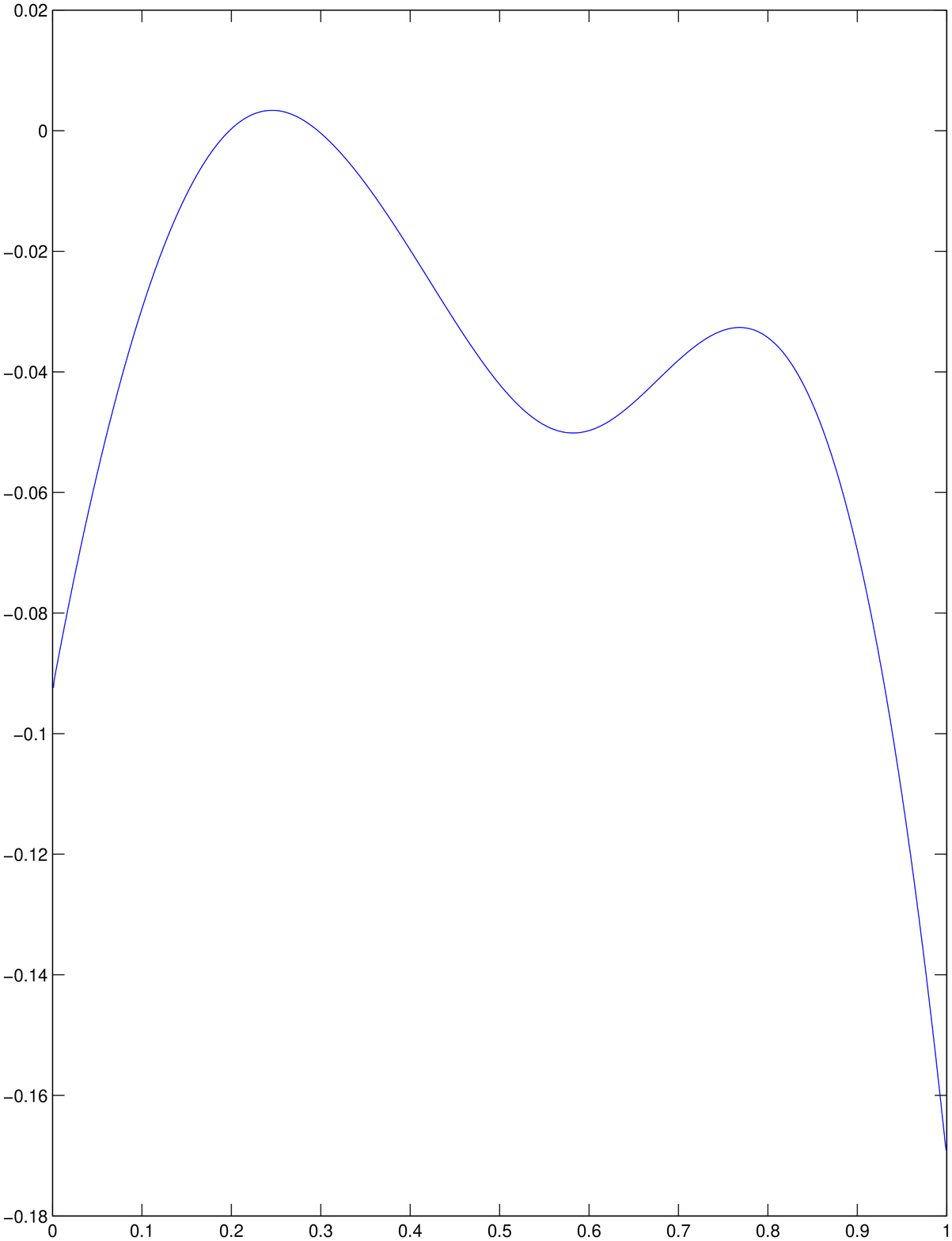} \; \includegraphics[height=45mm,width=50mm]{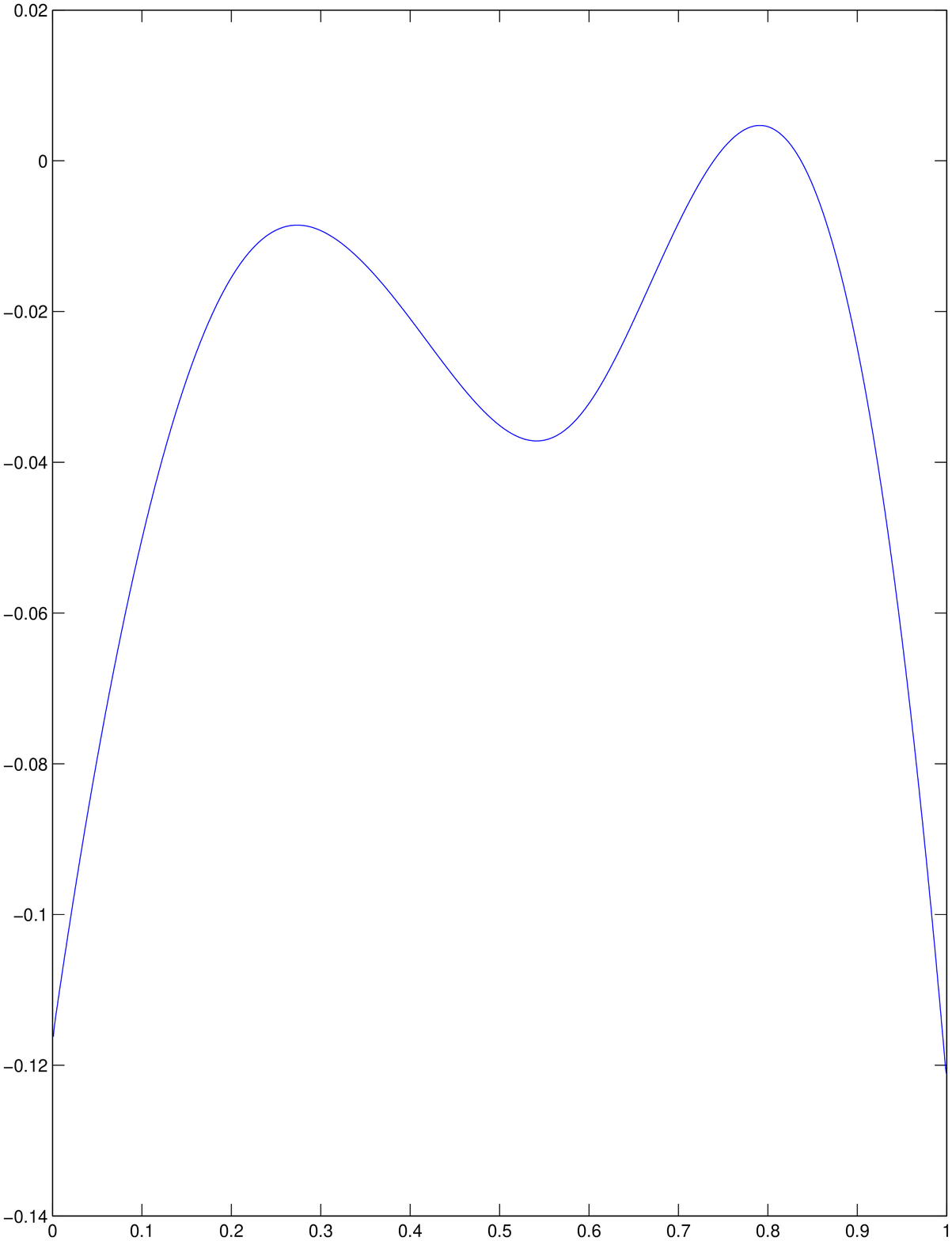}\; \includegraphics[height=45mm,width=50mm]{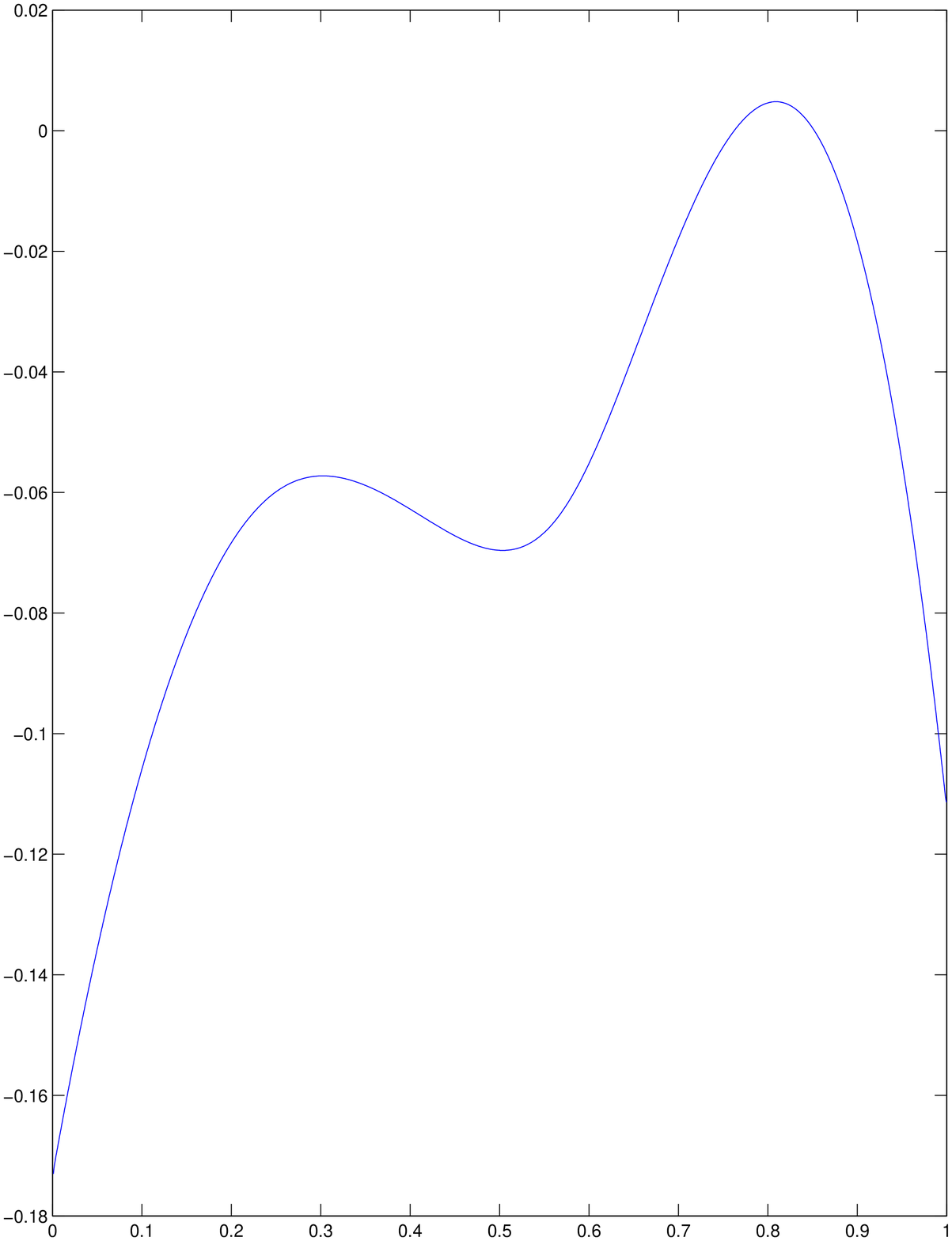} 
\end{center}
\vspace{-.7cm}
\caption{
\label{fig:discphi} 
{\footnotesize\sc 
Discontinuity in example \fer{exR} as in Figure \ref{fig:disc}. The phase $\vp(t,x)$ as a function of $x$ at different times around the discontinuity time. We observe that the discontinuities depicted in Figure \ref{fig:disc} correspond to a smooth transition on $\vp$. 
}
        }
\end{figure}

%-------------------------------------------------------------
\subsection{A numerical counter-example}
%-------------------------------------------------------------

We give a numerical simulation of a similar phenomena which however does not use a discontinuity of the gradient to generate the discontinuities of the Lagrange multiplier and the concentration point. We use again $\psi \equiv 1$ ($I= \varrho$), and 
\beq
R(x, \varrho) = (x -x^2 +3 x^4) (9-(1+x)^3)  -  \varrho. 
\label{exR}
\eeq
The computations depicted in Figures \ref{fig:disc}  and \ref{fig:discphi}  are performed with $0 \leq x \leq 1$ and $\vp^0(x)=-|x-0.05|$. They show how a smooth change on the phase funtion $\vp(t,x)$ along time, can produce a discontinuity of the Lagrange multiplier $I(t)$ and the concentration point $\bar x (t)$. A local maximum of $\vp(t,x)$ in $x$ becomes larger as time goes and becomes the global maximum at the discontinuity time. This follows the scenario in the analytical counter-example.

%-------------------------------------------------------------
\subsection{Concavity conditions for continuity}
%-------------------------------------------------------------

In this section, we are interested in conditions that imply the continuity of the $I(t)$ and $\bar x (t)$. Because a discontinuity corresponds to a double maximum of $\vp(t,x)$, it is natural to look for conditions implying the concavity of the solution $\vp$. This is the purpose of the next 
%-----------------------------------
\begin{theorem} Assume \fer{asr2} and that $R(\cdot,I)$ is strictly concave, uniformly for bounded $I$. 
Then, for any $\vp^0 \in W^{1, \infty}(\R^d)$, $\vp^0$ uniformly concave, any solution to \fer{eq:chj} is strictly concave and thus  $\bar x (t)$ is continuous. Therefore, in Theorem \ref{th:lscx}, the equalities hold everywhere and $I(t)$ is also continuous.
\label{th:continuity}
\end{theorem}

\proof Let $\vp$ a solution of \fer{eq:chj} and $0 < \alpha <1$. We consider the function
$$ \chi(x,y,t) := \alpha \vp (t,x) + (1-\alpha) \vp (t,y) - \vp (t,\alpha x + (1-\alpha)y)\; .$$
The aim is to prove that $\chi(x,y,t) < 0$ for any $x,y \in \R^d$ such that $x \neq y$ and for any $t$. To do so, we first show that $\chi$ is a viscosity subsolution of
$$ \chi_t \leq 2K |\chi_x + \chi_y| +  \f{\chi^2_x}{\alpha} +  \f{\chi^2_y}{1-\alpha} + \tilde R_\alpha (x,y,t) \quad \hbox{in  }\R^d \times \R^d \times (0, +\infty)\; ,
$$
where $K:=||D\vp ||_\infty$ and
$$ \tilde R_\alpha (x,y,t) = \max_{|J|Ê\leq ||I||_\infty} \left\{\alpha R (x, J) + (1-\alpha)R (y, J) - R (\alpha x + (1-\alpha) y, J)\right\}\; .$$
We are going to argue below as if the function $t \mapsto R \big(x, I(t)\big)$ were continuous : a rigourous proof would consists in approximating $I$ by continuous functions, in proving that the corresponding $\chi$'s satisfy the above inequality (this is also why we put a max in the definition of $ \tilde R_\alpha$ :  to point out the uniformity in $I$) and finally to pass to the limit. We drop these details here for the sake of simplicity in the presentation.

Consider a test function $\phi(x,y,t)$ in the viscosity sense and  $(\xb, \yb, \tb)$ is a strict maximum point of $\chi(x,y,t) - \phi (x,y,t)$, we look for local maximum points of
$$ \alpha \vp (t,x) + (1-\alpha) \vp (t,y) - \vp (t,z) - \f{|\alpha x + (1-\alpha)y-z|^2}{\e}- \phi (x,y,t) $$
near $(\xb, \yb, \tb)$. By classical arguments, since $(\xb, \yb, \tb)$ is a {\em strict} 
maximum point of $\chi(x,y,t) - \phi (x,y,t)$, there exists a sequence $(\xe, \ye, ,\ze, \te)$ of maximum points of this new function which converge to $(\xb, \yb, \alpha \xb + (1-\alpha)\yb, \tb)$ as $\e$ tends to $0$.

If $p :=  \f{2(\alpha x + (1-\alpha)y-z)}{\e}$, standard uniqueness arguments (cf. \cite{CIL}) provide the existence of $a,b,c$ (playing the roles of the $t$-derivatives of $\vp$ at respectively $(\xe,\te)$, $(\ye,\te)$, $(\ze,\te)$) such that
$$ a \leq |Êp + \f{D_x \phi (\xe,\ye, \te)}{\alpha}|^2 + R \big(\xe, I(\te)\big) \; ,$$
$$ b \leq |Êp + \f{D_y \phi (\xe,\ye, \te)}{1-\alpha}|^2 + R \big(\ye, I(\te)\big) \; ,$$
$$ c \geq |Êp |^2 + R \big(\ze, I(\te)\big) \; , $$
and with $\alpha a + (1-\alpha)b - c \geq \phi_t (\xe,\ye,\te)$.

Combining these inequalities, we obtain, after straightforward computations
$$ \phi_t (\xe,\ye,\te) \leq 2Êp \cdot (D_x \phi (\xe,\ye, \te) + D_y \phi (\xe,\ye, \te)) + 
$$ 
$$
 \left| \f{D_x \phi (\xe,\ye, \te)}{\alpha}\right|^2 + \left|\f{D_y \phi (\xe,\ye, \te)}{1-\alpha}\right|^2 + \alpha R \big(\xe, I(\te)\big) 
 + (1-\alpha)R \big(\ye, I(\te)\big) - R\big(\ze, I(\te)\big) \; .
 $$
To conclude, we first remark that $|p| \leq ||D\vp ||_\infty$ since $p$ is in the $x$-subdifferential of $\vp$ at $x$, and then let $\e$ tends to $0$. Since $\ze \to \alpha \xb + (1-\alpha) \yb$, we have 
$$\alpha R \big(\xe, I(\te)\big) + (1-\alpha)R \big(\ye, I(\te)\big) - R \big(\ze, I(\te)\big) \to \alpha R \big(\xb, I(\tb)\big) + (1-\alpha)R \big(\yb, I(\tb)\big) - R \big(\alpha \xb + (1-\alpha) \yb, I(\tb)\big) 
$$ and we are done.

In order to prove that $\chi(x,y,t) < 0$ for any $x,y \in \R^d$ such that $x \neq y$ and $t>0$, we first remark that $\chi(x,y,t) \leq  0$ for any $x,y \in \R^d$ as a consequence of a standard comparison result~: indeed since $\tilde R_\alpha (x,y,t) \leq 0 $ (because $R(\cdot,I)$ is concave for any $I$), $\chi$ is a subsolution of the P. D. E. 
$$ w_t = 2K |w_x + w_y| +  \f{w^2_x}{\alpha} +  \f{w^2_y}{1-\alpha} \quad \hbox{in  }\R^d \times \R^d \times (0, +\infty)\; ,$$
while $0$ is a solution, and $\chi(x,y,0) \leq  0$ in $ \R^d \times  \R^d$ because $\vp^0$ is strictly concave; therefore  $\chi(x,y,t) \leq  0$ in $\R^d \times \R^d \times (0,+\infty)$. 

On the other hand, if $\chi(x,y,t) = 0$ at some point $(x,y,t)$, then $(x,y,t)$ is a maximum point of $\chi-0$ and the viscosity subsolution inequality reads
$$ 0 \leq \tilde R_\alpha (x,y,t) \; ,$$
a contradiction with the strict concavity of $R(\cdot,I)$, uniformly in $I$ for $I$ bounded, which implies $\tilde R_\alpha (x,y,t) < 0$.
 \qed

%%%%%%%%%%%%%%%%%%%%%%%%%%%%%%%%%%%%%%%%%%%%
\section{The monomorphic case: more precise statements}
\label{sec:mono}
%-------------------------------------------
%%%%%%%%%%%%%%%%%%%%%%%%%%%%%%%%%%%%%%%%%%%%
%%%%%%%%%%%%%%%%%%%%%%%%%%%%%%%%%%%%%%%%%%%%

We consider in this section the special case when $d=1$ and we assume a {\em monomorphic} situation according to \fer{as:mono} which leads us to introduce the following assumption

There exists a function $\tilde R(x,I)$ such that $ x \mapsto \tilde R(x,I)$ is strictly monotone for any $I$, $ I \mapsto \tilde R(x,I)$ is strictly decreasing for any $x$ and
\begin{equation}\label{hrt}
\tilde R(x,I) = 0\  (resp.\ >0,\ <0)\quad \hbox{if and only if}  \quad R(x,I) = 0\  (resp.\ >0,\ <0).
\end{equation}

A typical example of such situation was given in the introduction, see \fer{mie} with assumptions \fer{asq1}--\fer{as:lip}. When $R$ is of the form $R(x,I) = b(x) Q_1(I) -d(x)  Q_2(I)$. Then \fer{hrt} reduces to 
\beq
x \mapsto \f{b(x)}{d(x)}  \text{ is strictly monotone}. 
\label{eq:mono}
\eeq
Indeed, one can choose $\displaystyle \tilde R(x,I) = \f{R(x,I)}{d(x) Q_1(I)} =  \f{b(x)}{d(x)} - \f{Q_2(I)}{Q_1(I)}$.
This is weaker than assuming $b' >0$ and $d'<0$ for which $R$ satifies  \fer{hrt}.

%---------------------------------------------------------------
\begin{theorem} Assume \fer{as:lip}, \fer{hrt} and that $\vp^0 \in W^{1, \infty}(\R^d)$ satisfies
\begin{equation}\label{infty}
 \max_{\R} \vp^0 (x)= 0 \quad \hbox{and} \quad \limsup_{|x| \to + \infty}\, \vp^0 (x) < 0\; .
\end{equation}
For all $t>0$ except perhaps for a countable number, there exists a unique $x(t)$ such that
$$
 \vp(t,\bar x(t))=\max_{\R} \vp (t,x)= 0\; .
$$
Moreover $t \mapsto x(t)$ is monotone and, with the notations  $x(t^+) = \lim_{s\downarrow t} x(s)$, $x(t^-) = \lim_{s \uparrow t} x(s)$, then
\begin{enumerate}
\item If $\tilde R$ is increasing in $x$ then $x(t^+)$ the largest maximum point of $\vp(t,\cdot)$ and $x(t^-)$ is the smallest one.
\item If $\tilde R$ is decreasing in $x$ then $x(t^+)$ the smallest maximum point of $\vp(t,\cdot)$ and $x(t^-)$ is the largest one.
\end{enumerate}
\end{theorem}
%----------------------------------------------------------------

\proof Using assumption (\ref{infty}) and cone of dependence type properties (recall that $\vp$ is Lipschitz continuous), we see that, for any $T>0$ and $t \in (0,T)$, the maximum of the function $x\mapsto \vp (t,x)$ is achieved in a fixed compact subset of $\R$. And, by (\ref{eq:iv}), if $t$ is a Lebesgue's point of $I$ (meaning here a continuity point of the increasing function $I$) and $\bar x$ is maximum of $\vp (\cdot,t)$, we have
$$ R\big(\bar x, I (t)\big) = 0\; .$$
Therefore, by \fer{hrt},
$$ \tilde R\big(\bar x, I (t)\big) = 0\; ;$$
we deduce that $\bar x$ is unique for such $t$, we denote it by $x(t)$ and, from (\ref{asq1}) and (\ref{eq:mono}), we see that $t \mapsto x(t)$ has the same monotonicity as $\tilde R$ since $I(t)$ is increasing.

If $t_k \downarrow t$, passing to the limit in $\vp(t_k, x(t_k)) = 0$, we see that $\vp(t, x(t^+)) = 0$. Choosing a sequence $(t_k)_k$ of continuity points of $x(\cdot)$ or equivalently of $I(\cdot)$, we have also
$$ \tilde R\big(x(t^+), I(t^+)\big) = 0 \; .$$

Now assume that $\tilde R$ is strictly increasing in $x$ and that there exists a maximum point $\bar x > x(t^+)$ of $\vp(t,\cdot)$. Then for $y$ close enough to $\bar x$ and $s>t$, using that $I (t^+) = \lim_{s\downarrow t}\, I(s)$, we have
$$ \tilde R \big(y,I (s)\big) \geq \tilde R \big(x(t^+),I (s)\big)+\eta\; ,$$
for some $\eta >0$, and therefore, using again \fer{hrt}, for $s$ close enough to $t$
$$ R \big(y,I (s)\big) \geq \tilde \eta\; ,$$
for some $\tilde \eta$, which would imply that $\vp(s,\bar x)>0$, a contradiction. Therefore such $\bar x$ does not exist and actually, $x(t^+)$ is the largest maximum point of $\vp(t,\cdot)$.

In the same way, $x(t^-)$ is the smallest maximum point of $\vp(t,\cdot)$ and the case when $\tilde R$ is strictly decreasing follows along the same lines.

%%%%%%%%%%%%%%%%%%%%%%%%%%%%%%%%%%%%%%%%%%%%
\section{Several environmental unknowns $I^{(k)}_\e(t)$}
\label{sec:seu}
%-------------------------------------------
%%%%%%%%%%%%%%%%%%%%%%%%%%%%%%%%%%%%%%%%%%%%

Several possible extensions are possible which are more realistic and give more interesting structures. In particular the selection rate $R$ may depend on several integrals and competition may be modeled by convolution terms (\cite{genieysP, DJMR, MG}). Here we study the former extension when the model contains several environmental unknowns. This means that we consider the equation
\beq \left\{ \begin{array}{l}
\f{\p}{\p t} n_\e - \e \Delta n_\e =\f{n_\e}{\e} R\big(x,I_\e(t)\big), \qquad x\in \R^d, \; t \geq 0,
\\ \\
n_\e(t=0)= n^0_\e \in L^1(\R^d) , \quad n^0_\e\geq 0,
\end{array} \right.
\label{eq:masterext}
\eeq
and now the environment is described by a vector valued unknown $I_\e= (I_\e^{(1)},..., I_\e^{(I)})$, with 
\beq 
I^{(i)}_\e(t)= \int \psi^{(i)}(x) n_\e(t,x) dx, \qquad i=1,...,N.
\label{eq:envirext}
\eeq
The  $\psi^{(i)}(x)$ are still given positive functions (see the precise assumptions below).
\\

We may have in mind the example when mortality and birth terms depend differently on the total population, i.e.,    $N=2$ and 
$$
R(x, I)= b(x) Q_1\big(I^{(1)}(t)\big) -d(x)  Q_2 \big(I^{(2)}(t)\big).
$$
We may also have in mind the example of the chemostat with several nutrients
$$
R(x, I)= \sum_{i=1}^I b_i(x) Q_i \big(I^{(i)}(t)\big) -d(x).
$$

Our first goal is to derive the $L^\infty$ estimates. This uses the following assumptions; there are constants $\psi_m<\psi_M$, $0< I_m \leq  I_M < \infty$ such that (vector inequalities should be understood componentwise)
\beq
0<\psi_m \leq \psi \leq \psi_M < \infty,  \qquad \psi^{(i)} \in W^{2,\infty}(\R^d) . 
\label{as:sev1}
\eeq
\beq
\min_{x\in\R^d} R(x, I) \leq 0, \; \text{ when  } \; \min_{1\leq i \leq N } I^{(i)} \geq I_M,
\label{as:sev2}
\eeq
\beq
\max_{x\in\R^d} R(x, I) \geq 0, \; \text{ when  } \;  \max_{1\leq i \leq N } I^{(i)}  \leq I_m,
\label{as:sev3}
\eeq
and there exists a constant $K>0$ such that, for any $x\in \R^d,$ $I \in \R^N$
\beq
-K \leq \f{\p}{\p I^{(i)}} R(x,I) < - K^{-1} <0,    \qquad \sup_{  I_m/2 \leq I \leq 2 I_M} \|R(\cdot, I)\|_{W^{2,\infty}(\R^d)} \leq K.  
\label{as:sev4}
\eeq

Then, an easy adaptation of the proof of Theorem \ref{th:basic} yields
%-----------------------------
\begin{theorem}\label{BdSeveral} With the assumptions \fer{as:sev1}--\fer{as:sev4} and $I_m - C \e^2\leq I_\e(0) \leq I_M+ C \e^2$ (componentwise), we have the uniform bound 
$$
I_m - C \e^2 \leq I^{(i)}(t) \leq I_M +C \e^2,  \qquad \forall t \geq 0, \quad \forall i=1,...,N,
$$
with $C$ a  positive constant. 
\end{theorem}

It is unclear if one can expect $BV$ estimates (and monotonicity results in the limit) on the quantities $I^{(i)}_\e$.
Assuming such a bound, one might follow the arguments of Section \ref{sec:aa} and derive again the constrained H.-J. equation \fer{eq:chj} where now $I(t)$ is a $N$ dimensional Lagrange multiplier associated with the (single) constraint $\max_{x\in \R^d} \vp(t,x) =0$. Of course this rises several fundamental questions concerning uniqueness of the solution to this constrained H.-J. equation. Some information is clearly lost because the quantities $I^{(i)}$ are related together in that case when  $n_\e$ concentrates on less than $N$ Dirac masses (for instance in the monomorphic case).  Therefore we face here a fundamental non-uniqueness situation.

\end{document}